\documentstyle[11pt]{article}

\newfont{\tenmsb}{msbm10 scaled\magstep1}

\let\ssection=\section\renewcommand{\section}{\setcounter{equation}{0}\ssection}

\def\smallover#1/#2{\hbox{$\textstyle{#1\over#2}$}}

\begin{document}

\title{A note on the generalized Weierstrass representation}
\author{L.Martina$^{1}$, Kur.Myrzakul$^{2}$\thanks{%
Permanent address: Institute of Mathematics, National Academy of Sciences,
Alma-Ata, Kazakhstan}, R.Myrzakulov$^{2}$ \\
$^{1}$Dipartimento di Fisica dell'Universit\`{a} and Sezione INFN di Lecce,\\
via Arnesano,73100 Lecce, Italy\\
$^{2}$Institute of Physics and Technology, 480082 Alma-Ata-82, Kazakhstan}
\date{\today}
\maketitle

\begin{abstract}
The study of the relation between the Weierstrass inducing formulae for
constant mean curvature surfaces and the completely integrable euclidean
nonlinear $\sigma $-model suggests a connection among integrable $\sigma $%
-models in a background and other type of surfaces. We show how a
generalization of the Weierstrass representation can be achieved and we
establish a connection with the Weingarten surfaces. We suggest also a
possible generalization for two-dimensional surfaces immersed in a flat
space $R^{8}$ with Euclidean metric.
\end{abstract}


\section{ Introduction}

It is well known that conformal immersions of the minimal surfaces into a
3-dimensional euclidean space are described by a classical Weierstrass
representation (WR) 
$$
x^{+}=2i\int_{z_{0}}^{z}(\bar{\psi}_{1}^{2}dz^{\prime }-\bar{\psi}_{2}^{2}d%
\bar{z}^{\prime }),\quad x^{-}=2i\int_{z_{0}}^{z}(\psi _{2}^{2}dz^{\prime
}-\psi _{1}^{2}d\bar{z}^{\prime }),\quad x_{3}=-2\int_{z_{0}}^{z}(\psi _{2}%
\bar{\psi}_{1}dz^{\prime }+\psi _{1}\bar{\psi}_{2}d\bar{z}^{\prime })\eqno(1)
$$
where $x^{\pm }=x_{1}\pm ix_{2}$ and $\bar{\psi}_{1},$ $\psi _{2}$ are
arbitrary analytic functions. Generalization of WR was proposed by
Konopelchenko in [1] (see, also [2-4]) in order to study constant mean
curvature surfaces. In this case $\psi _{1}$, $\psi _{2}$ satisfy the
following infinite-dimensional Hamiltonian system 
$$
\psi _{1z}=2H(|\psi _{1}|^{2}+|\psi _{2}|^{2})\psi _{2},\qquad \psi _{2\bar{z%
}}=-2H(|\psi _{1}|^{2}+|\psi _{2}|^{2})\psi _{1}.\eqno(2) 
$$
One can assume, without loss of generality $H=\frac{1}{2}$. Then, system (1)
takes the 2-dimensional Dirac-like equation 
$$
\psi _{1z}=u\psi _{2},\qquad \psi _{2\bar{z}}=-u\psi _{1}\eqno(3) 
$$
where 
$$
u=|\psi _{1}|^{2}+|\psi _{2}|^{2}.\eqno(4) 
$$
Using the standard formulas, we find that the first fundamental form on the
surface is given by 
$$
I=u^{2}dzd\bar{z}\eqno(5) 
$$
at this the Gassian curvature is 
$$
K=-\frac{4}{u^{2}}[logu]_{z\bar{z}}.\eqno(6) 
$$

In this paper we would like to consider a generalization of the equation
(3),  study its connection with certain nonlinear $\sigma $ - models and the
corresponding surfaces in a 3-dimensional Euclidean space, obtained via a
suitable Weierstrass-like formula.

\section{Generalized Weierstrass representation}

First, differently from the system (3), we postulate the equations
$$
\psi _{1z}=u\psi _{2}+h_{1}\psi _{1},\qquad \psi _{2\bar{z}}=-u\psi
_{1}+h_{2}\psi _{2},\eqno(7a)
$$
$$
\bar{\psi}_{1\bar{z}}=u\bar{\psi}_{2}+\bar{h}_{1}\bar{\psi}_{1},\qquad \bar{%
\psi}_{2z}=-u\bar{\psi}_{1}+\bar{h}_{2}\bar{\psi}_{2},\eqno(7b)
$$
where $h_{j}$ are complex functions to be characterized by further
requirements. We call the equations (7) the generalized Konopelchenko system
(GKS). Second, we postulate that the generalized Weierstrass representation
(GWR) has the form 
$$
x^{+}=2\int_{\gamma }F^{+}(z,\bar{z})dz+G^{+}(z,\bar{z})d\bar{z}\eqno(8a)
$$
$$
x^{-}=2\int_{\gamma }F^{-}(z,\bar{z})dz+G^{-}(z,\bar{z})d\bar{z}\eqno(8b)
$$
$$
x_{3}=-2\int_{\gamma }F_{3}(z,\bar{z})dz+G_{3}(z,\bar{z})d\bar{z}\eqno(8c)
$$
where 
$$
F^{+}=\rho (\bar{\psi}_{1}^{2}-g\bar{\psi}_{2}^{2}),\quad F^{-}=\rho (\psi
_{2}^{2}-g\psi _{1}^{2}),\quad F_{3}=\rho (\bar{\psi}_{1}\psi _{2}+g\psi _{1}%
\bar{\psi}_{2})\eqno(9a)
$$
$$
G^{+}=\bar{\rho}(\bar{\psi}_{2}^{2}-\bar{g}\bar{\psi}_{1}^{2}),\quad G^{-}=%
\bar{\rho}(\psi _{1}^{2}-\bar{g}\psi _{2}^{2}),\quad G_{3}=\bar{\rho}(\psi
_{1}\bar{\psi}_{2}+\bar{g}\bar{\psi}_{1}\psi _{2})\eqno(9b)
$$
with $\rho ,g$ are complex functions to be defined and $\gamma $ is any
curve from a fixed point to $z$. Since the coordinates $x_{i}$ do not depend
on the choice of the curve $\gamma $ in the complex plane $C$ (but only its
end points), the functions $F_{i},G_{i}$ must satisfy the following
conditions 
$$
F_{\bar{z}}^{+}=G_{z}^{+},\quad F_{\bar{z}}^{-}=G_{z}^{-},\quad F_{3\bar{z}%
}=G_{3z}.\eqno(10)
$$
These show that the integrands are total differentials. Substituting (9)
into (10), we get the following restrictions for the functions $h_{i},\rho $
and $g$ 
$$
\lbrack \rho _{\bar{z}}+2\rho \bar{h}_{1}+(\bar{\rho}\bar{g})_{z}]\bar{\psi}%
_{1}^{2}-[\bar{\rho}_{z}+2\bar{\rho}\bar{h}_{2}+(\rho g)_{\bar{z}}]\bar{\psi}%
_{2}^{2}+2u(\rho +\bar{\rho})\bar{\psi}_{1}\bar{\psi}_{2}+2(\bar{\rho}\bar{g}%
\bar{\psi}_{1}\bar{\psi}_{1z}-\rho g\bar{\psi}_{2}\bar{\psi}_{2\bar{z}})=0,%
\eqno(11a)
$$
$$
\lbrack \rho _{\bar{z}}+2\rho h_{2}+(\bar{\rho}\bar{g})_{z}]\psi _{2}^{2}-[%
\bar{\rho}_{z}+2\bar{\rho}h_{1}+(\rho g)_{\bar{z}}]\psi _{1}^{2}-2u(\rho +%
\bar{\rho})\psi _{1}\psi _{2}-2(\rho g\psi _{1}\psi _{1\bar{z}}-\bar{\rho}%
\bar{g}\psi _{2}\psi _{2z})=0,\eqno(11b)
$$
$$
\lbrack \rho _{\bar{z}}+\rho (\bar{h}_{1}+h_{2})]\bar{\psi}_{1}\psi _{2}-[%
\bar{\rho}_{z}+\bar{\rho}(h_{1}+\bar{h}_{2})]\psi _{1}\bar{\psi}_{2}+u(\rho -%
\bar{\rho})(|\psi _{2}|^{2}-|\psi _{1}|^{2})+(\rho g\psi _{1}\bar{\psi}%
_{2})_{\bar{z}}-(\bar{\rho}\bar{g}\bar{\psi}_{1}\psi _{2})_{z}=0.\eqno(11c)
$$

\section{The associated generalized $\sigma$-model of the GKS}

In the paper [6] it was shown that the Konopelchenko system (3) is
equivalent to a scalar second order equation, which is the stereographic
projection of the euclidean nonlinear $\sigma $ - model. Inspired by the
transformation used there, also in this we define a function $\omega $

$$
\omega =\frac{\psi _{1}}{\bar{\psi}_{2}}.\eqno(12) 
$$

Now, computing the derivatives of $\omega $ with the use of system (7) we
find 
$$
\omega _{z}=\psi _{2}^{2}\left( 1+\left| \omega \right| ^{2}\right) +h\omega
,\eqno(13)
$$
where we have introduced the notation $h=h_{1}-\bar{h}_{2}.$ Solving the
previous relation with respect to $\psi _{2}^{2}$ , we obtain an ''inverse''
formulas for $\psi _{1},\psi _{2},$ i.e.

$$
\psi _{1}=\epsilon \omega \frac{\sqrt{\bar{\omega}_{\bar{z}}-\bar{h}\bar{%
\omega}}}{1+|\omega |^{2}},\qquad \psi _{2}=\epsilon \frac{\sqrt{\omega
_{z}-h\omega }}{1+|\omega |^{2}},\quad \epsilon ^{2}=1.\eqno(14)
$$
The formulae (14) generalize those used in [6] and they coincide with for $%
h=0.$ Then, system (7) allows to find an equation for $\omega $,  precisely  
$$
\omega _{z\bar{z}}=2\frac{\omega _{z}\omega _{\bar{z}}}{1+|\omega |^{2}}\bar{%
\omega}+\left( h\omega \right) _{\bar{z}}+2\left( \omega _{z}-h\omega
\right) h_{2}-2\frac{\left| \omega \right| ^{2}}{1+|\omega |^{2}}\left[
\omega _{\bar{z}}h-\omega _{z}\bar{h}+\omega \left| h\right| ^{2}\right] .%
\eqno(15)
$$
We can show that in the  case $h=0$ the equation (15) fits both the usual
euclidean nonlinear $\sigma$ - model  and the Ernst-like system
studied by Schief in [11] .

\section{The $O(3)$ $\sigma $-model}

Now, let us consider the simplest reductions of the above equations.

i) Let us consider first the case $h_{2}=\bar{h}_{1}=g=0,\quad \rho =i$.
Then, the conditions (11) are obviously satisfied and the GKS (7) reduces to
the usual KS (3). Moreover, the system (15) reduces to the usual euclidean $%
O(3)$ $\sigma $-model 
$$
\omega _{z\bar{z}}=2\frac{\omega _{z}\omega _{\bar{z}}}{1+|\omega |^{2}}\bar{%
\omega}\eqno(16)
$$
the integrability properties of which are well known [12] and its connection
connected with the constant mean curvature surfaces is described in  [6],
[7] and [9].

ii) Now let $h_{2}=\bar{h}_{1}=0,\quad \rho =i$. The restrictions (11) now
dictate strong conditions on the function $g$, i.e. 
$$
(\bar{g}\bar{\psi}_{1}^{2})_{z}+(g\bar{\psi}_{2}^{2})_{\bar{z}}=0,\eqno(17a)
$$
$$
(\bar{g}\psi _{2}^{2})_{z}+(g\psi _{1}^{2})_{\bar{z}}=0,\eqno(17b)
$$
$$
(g\psi _{1}\bar{\psi}_{2})_{\bar{z}}+(\bar{g}\bar{\psi}_{1}\psi _{2})_{z}=0.%
\eqno(17c)
$$
This equations were found in [9] and they have as a solution $g=J(z)/u^{2}$,
where $J(z)$ is an arbitrary analytic function. The corresponding
Weierstrass formula provides a non isothermic set of coordinates for the
mean constant curvature surfaces, since the corresponding scalar equation is
still (15).

\section{The $\sigma$-model in a curved space}

Now we would like to consider the restriction 
$$
h_{2}=\bar{h}_{1}=-\frac{\rho _{\bar{z}}\omega _{z}+\rho _{z}\omega _{\bar{z}%
}}{4\rm Re}(\rho )\omega _{z},\eqno(18)
$$
where the real function $\rho $ satisfies the equation 
$$
\rho _{z\bar{z}}=0.\eqno(19)
$$
In this case we have $h=0$ and the equation (15) becomes 
$$
\omega _{z\bar{z}}=2\frac{\omega _{z}\omega _{\bar{z}}}{1+|\omega |^{2}}\bar{%
\omega}+2h_{2}\omega _{z}.\eqno(20)
$$
Furthermore, the conditions (11) provide a constraint on $g$. Equation (20) 
represents  a nonlinear $\sigma $-model in a curved space, and it can be
interpreted also as a 2-dimensional isotropic magnet, whose squared
saturation moment is a harmonic function [5]. Moreover, it is known that
this system is integrable [5]. In fact, the equation (20) is the
compatibility condition for the linear $2\times 2$ matrix system

$$
\Psi _{z}=U\Psi ,\qquad \Psi _{z}=V\Psi \eqno(21) 
$$
where the potential matrices are given by

$$
U=\frac{\rho }{\varrho +\rho }S_{z}S,\qquad V=-\frac{\rho }{\varrho -\rho }%
S_{z}S\eqno(22)
$$
with 
$$
S=\frac{i}{1+|\omega |^{2}}\left( 
\begin{array}{ll}
1-|\omega |^{2} & -2\bar{\omega} \\ 
-2\omega  & |\omega |^{2}-1
\end{array}
\right) ,\quad \varrho =i\beta -u+\sqrt{(u-\gamma )(u+\gamma )},\qquad
\gamma =\rho +i\beta .\eqno(23)
$$
The function $\gamma $ is analytic in a suitable region and $u\in C$ is a
''hidden'' spectral parameter, which  does not depend on the coordinates.

\section{The Ernst-type equation}

Weakening the reality condition on $\rho $ used in the previous section, we
consider the case 
$$
h_{2}=\bar{h}_{1}=-\frac{p_{\bar{z}}\omega _{z}+p_{z}\omega _{\bar{z}}}{4%
{\rm Re}(p)\omega _{z}}=\frac{f}{2\omega _{z}},\qquad p_{z\bar{z}}=0,%
\eqno(24)
$$
where we introduced $p=\rho +i\sigma $ , both these functions are harmonic,
and 
$$
f=-\frac{1}{2}\frac{p_{\bar{z}}\omega _{z}+p_{z}\omega _{\bar{z}}}{{\it Re}%
(p)}.\eqno(25)
$$
Now the equation (15) takes the form 
$$
\omega _{z\bar{z}}=2\frac{\omega _{z}\omega _{\bar{z}}}{1+|\omega |^{2}}\bar{%
\omega}+f,\qquad p_{z\bar{z}}=0,\eqno(26)
$$
and it was considered by Schief [11] as a particular reduction of the static
Loewner system in connection with the study of the Weingarten surfaces of
Class 2. The system is integrable and this property was studied in [13].
With these assumptions the GKS  (7) takes the form 
$$
\psi _{1z}=u\psi _{2}+\frac{\bar{f}}{2\bar{\omega}_{\bar{z}}}\psi
_{1},\qquad \psi _{2\bar{z}}=-u\psi _{1}+\frac{f}{2\omega _{z}}\psi _{2}%
\eqno(27a)
$$
$$
\bar{\psi}_{1\bar{z}}=u\bar{\psi}_{2}+\frac{f}{2\omega _{z}}\bar{\psi}%
_{1},\qquad \bar{\psi}_{2z}=-u\bar{\psi}_{1}+\frac{\bar{f}}{2\bar{\omega}_{%
\bar{z}}}\bar{\psi}_{2}\eqno(27b)
$$
where 
$$
u=|\psi _{1}|^{2}+|\psi _{2}|^{2}=\frac{(\omega _{z}\bar{\omega}_{\bar{z}%
})^{1/2}}{1+|\omega |^{2}}.\eqno(28)
$$

We notice that both the systems (26) and (27) are conformally invariant
under the transformations 
\begin{eqnarray*}
z^{\prime } &=&\xi \left( z\right) ,\qquad \bar{z}^{\prime }=\bar{\xi}\left( 
\bar{z}\right) , \\
\omega ^{\prime } &=&\omega ,\qquad p^{\prime }=p,\qquad \psi _{1}^{\prime
}=\left( \bar{\xi}_{\bar{z}}\right) ^{-\frac{1}{2}}\psi _{1}^{{}},\qquad
\psi _{2}^{\prime }=\left( \xi _{z}\right) ^{-\frac{1}{2}}\psi _{2}^{{}},
\end{eqnarray*}
where $\xi \left( z\right) $ is an arbitrary analytic function on the
complex plane. This property is important from the geometrical point of
view, since it is the realization of the invariance under parametrization of
the 2-dimensional surfaces, which are in the motivations of the present
work. Furthermore, this fact indicates that the system admits conservation
laws, which would be useful in the further analysis.

Now, we consider certain reductions of the system (27). Besides the trivial
example $p=const$, which leads to the usual system (3), we can consider the
following special cases.

i) Case $p=p(\bar{z})$ leads to the quasi-linear system with peculiar
non-constant coefficients 
$$
\psi _{1z}=u\psi _{2}-\frac{\bar{p}_{z}}{4Re(p)}\psi _{1},\quad \psi _{2\bar{%
z}}=-u\psi _{1}-\frac{p_{\bar{z}}}{4Re(p)}\psi _{2}.\eqno(29) 
$$

ii) Case $p=p(z)$. In this case from (7) we obtain a system which is highly
nonlinear in the $\psi ^{\prime }s,$ and which can be written in form 
$$
\psi _{1z}=u\psi _{2}-\epsilon \omega \frac{\bar{p}_{\bar{z}}(\bar{\omega}_{%
\bar{z}})^{-\frac{1}{2}}\bar{\omega}_{z}}{4Re(p)(1+|\omega |^{2})},\quad
\psi _{2\bar{z}}=-u\psi _{1}-\epsilon \frac{p_{z}(\omega _{z})^{-\frac{1}{2}%
}\omega _{\bar{z}}}{4Re(p)(1+|\omega |^{2})}.\eqno(30)
$$

Our aim is now to deduce directly from the above expressions the integrals
(8), solving implicitly the system of constraints  (11). This result can be
achieved looking at a conservation law form of the system (26). In doing so,
we introduce the ''spin'' $S$ matrix in (23). Then, equation (26) is
equivalent, in the sense of the stereographic projection, to the matrix
equation

$$
(\rho [S,S_{\bar{z}}]+2\sigma _{\bar{z}}S)_{z}+(\rho [S,S_{z}]+2\sigma
_{z}S)_{\bar{z}}=0,\eqno(31)
$$
which has the conservation law form. For brevity we introduce the matrices 
$$
K(\rho ,\sigma ,\bar{z})=\rho [S,S_{\bar{z}}]+2\sigma _{\bar{z}}S,\quad
M(\rho ,\sigma ,z)=\rho [S,S_{z}]+2\sigma _{z}S,\eqno(32)
$$
which allows us to rewrite (31) in the form 
$$
K_{z}+M_{\bar{z}}=0.\eqno(33)
$$
The explicit expressions  of matrices $K$ and $M$ in terms of $\omega $ are 
$$
K=m_{{}}\left( 
\begin{array}{cc}
\bar{\omega}\omega _{\bar{z}}-\bar{\omega}_{\bar{z}}\omega  & -(\bar{\omega}%
_{\bar{z}}+\bar{\omega}^{2}\omega _{\bar{z}}) \\ 
\omega _{\bar{z}}+\omega ^{2}\bar{\omega}_{\bar{z}} & -(\bar{\omega}\omega _{%
\bar{z}}-\bar{\omega}_{\bar{z}}\omega )
\end{array}
\right) +n_{1}\left( 
\begin{array}{cc}
1-|\omega |^{2} & -2\bar{\omega} \\ 
-2\omega  & |\omega |^{2}-1
\end{array}
\right) \eqno(34a)
$$
and 
$$
M=m_{{}}\left( 
\begin{array}{cc}
\bar{\omega}\omega _{z}-\bar{\omega}_{z}\omega  & -(\bar{\omega}_{z}+\bar{%
\omega}^{2}\omega _{z}) \\ 
\omega _{z}+\omega ^{2}\bar{\omega}_{z} & -(\bar{\omega}\omega _{z}-\bar{%
\omega}_{z}\omega )
\end{array}
\right) +n_{2}\left( 
\begin{array}{cc}
1-|\omega |^{2} & -2\bar{\omega} \\ 
-2\omega  & |\omega |^{2}-1
\end{array}
\right) ,\eqno(34b)
$$
where 
$$
m_{{}}=-\frac{4\rho }{(1+|\omega |^{2})^{2}},\quad n_{1}=\frac{2i\sigma _{%
\bar{z}}}{1+|\omega |^{2}},\quad n_{2}=\frac{2i\sigma _{z}}{1+|\omega |^{2}}.%
\eqno(35)
$$
From the above expression it is easy to prove the properties 
$$
K(\bar{\rho},\bar{\sigma},\bar{z})=-M^{\dagger }(\rho ,\sigma ,z),\quad
K(\rho ,\sigma ,\bar{z})=-M^{\dagger }(\bar{\rho},\bar{\sigma}%
,z)=M^{{}}(\rho ,\sigma ,\bar{z}).\eqno(36)
$$
We need to use also the following relations 
$$
\omega _{z}=\frac{u^{2}}{\bar{\psi}_{2}^{2}},\qquad \bar{\omega}_{z}=-\frac{T%
}{\psi _{2}^{2}},\qquad T=\bar{\psi}_{1}\psi _{2z}-\bar{\psi}_{1z}\psi _{2}=2%
\frac{\omega _{z}\bar{\omega}_{z}}{(1+|\omega |^{2})^{2}}.\eqno(37)
$$
We note that $T$ satisfies the relation 
$$
T_{\bar{z}}=\frac{f}{\omega _{z}}T-\frac{u^{2}\bar{f}}{\bar{\omega}_{\bar{z}}%
},\eqno(38)
$$
which breaks the analiticity property possessed by the analogous quantity in
the theory of the constant mean curvature surfaces reviewed above [6, 9].
Using the formula (14) with $h=0$, we now find the expressions of matrices $K
$ and $M$ in terms of $\psi ^{\prime }$s 
$$
K=-4\rho \left( 
\begin{array}{cc}
-\psi _{1}\bar{\psi}_{2} & -\bar{\psi}_{2}^{2} \\ 
\psi _{1}^{2} & \psi _{1}\bar{\psi}_{2}
\end{array}
\right) +\frac{2i\sigma _{\bar{z}}}{u}\left( 
\begin{array}{cc}
|\psi _{2}|^{2}-|\psi _{1}|^{2} & -2\bar{\psi}_{1}\bar{\psi}_{2} \\ 
-2\psi _{1}\psi _{2} & |\psi _{1}|^{2}-|\psi _{2}|^{2}
\end{array}
\right) +\frac{4\rho \bar{T}}{u^{2}}\left( 
\begin{array}{cc}
\bar{\psi}_{1}\psi _{2} & -\bar{\psi}_{1}^{2} \\ 
\psi _{2}^{2} & -\bar{\psi}_{1}\psi _{2}
\end{array}
\right) \eqno(39)
$$
and 
$$
M=-4\rho \left( 
\begin{array}{cc}
\bar{\psi}_{1}\psi _{2} & -\bar{\psi}_{1}^{2} \\ 
\psi _{2}^{2} & -\bar{\psi}_{1}\psi _{2}
\end{array}
\right) +\frac{2i\sigma _{z}}{u}\left( 
\begin{array}{cc}
|\psi _{2}|^{2}-|\psi _{1}|^{2} & -2\bar{\psi}_{1}\bar{\psi}_{2} \\ 
-2\psi _{1}\psi _{2} & |\psi _{1}|^{2}-|\psi _{2}|^{2}
\end{array}
\right) +\frac{4\rho T}{u^{2}}\left( 
\begin{array}{cc}
-\psi _{1}\bar{\psi}_{2} & -\bar{\psi}_{2}^{2} \\ 
\psi _{1}^{2} & \bar{\psi}_{2}\psi _{1}
\end{array}
\right) .\eqno(40)
$$
In this notation the three conservation laws read 
$$
-[4\rho (\psi _{1}\bar{\psi}_{2}+\frac{\bar{T}}{u^{2}}\bar{\psi}_{1}\psi
_{2})+\frac{2i\sigma _{\bar{z}}}{u}(|\psi _{2}|^{2}-|\psi
_{1}|^{2})]_{z}+[4\rho (\bar{\psi}_{1}\psi _{2}+\frac{4\bar{T}}{u^{2}}\psi
_{1}\bar{\psi}_{2})-\frac{2i\sigma _{z}}{u}(|\psi _{2}|^{2}-|\psi
_{1}|^{2})]_{\bar{z}}=0,\eqno(41a)
$$

$$
\lbrack 4\rho (-\psi _{1}^{2}+\frac{\bar{T}}{u^{2}}\psi _{2}^{2})-\frac{%
4i\sigma _{\bar{z}}}{u}\psi _{1}\psi _{2}]_{z}+[4\rho (-\psi _{2}^{2}+\frac{%
4T}{u^{2}}\psi _{1}^{2})-\frac{4i\sigma _{z}}{u}\psi _{1}\psi _{2}]_{\bar{z}%
}=0,\eqno(41b)
$$
$$
\lbrack 4\rho \bar{\psi}_{2}^{2}-\frac{\bar{T}}{u^{2}}\bar{\psi}_{1}^{2})-%
\frac{4i\sigma _{\bar{z}}}{u}\bar{\psi}_{1}\bar{\psi}_{2}]_{z}+[4\rho (\bar{%
\psi}_{1}^{2}-\frac{4T}{u^{2}}\bar{\psi}_{2}^{2}-\frac{4i\sigma _{z}}{u}\bar{%
\psi}_{1}\bar{\psi}_{2}]_{\bar{z}}=0.\eqno(41c)
$$

As a result of these equations, we introduce the six real-valued functions 
$$
F_{1}=4\rho (\bar{\psi}_{1}^{2}-\frac{4T}{u^{2}}\bar{\psi}_{2}^{2})-\frac{%
4i\sigma _{z}}{u}\bar{\psi}_{1}\bar{\psi}_{2}+4\bar{\rho}(\psi _{2}^{2}-%
\frac{T}{u^{2}}\psi _{1}^{2})-\frac{4i\bar{\sigma}_{z}}{u}\psi _{1}\psi _{2},%
\eqno(42a)
$$
$$
F_{2}=4\rho (\bar{\psi}_{1}^{2}-\frac{4T}{u^{2}}\bar{\psi}_{2}^{2})-\frac{%
4i\sigma _{z}}{u}\bar{\psi}_{1}\bar{\psi}_{2}-4\bar{\rho}(\psi _{2}^{2}+%
\frac{T}{u^{2}}\psi _{1}^{2})+\frac{4i\bar{\sigma}_{z}}{u}\psi _{1}\psi _{2},%
\eqno(42b)
$$
$$
F_{3}=-2[4\rho (\bar{\psi}_{1}\psi _{2}+\frac{4\bar{T}}{u^{2}}\psi _{1}\bar{%
\psi}_{2})-\frac{2i\sigma _{z}}{u}(|\psi _{2}|^{2}-|\psi _{1}|^{2})],%
\eqno(42c)
$$
and
$$
G_{1}=-[4\bar{\rho}(\psi _{1}^{2}-\frac{4\bar{T}}{u^{2}}\psi _{2}^{2}-\frac{%
4i\bar{\sigma}_{\bar{z}}}{u}\psi _{1}\psi _{2}+4\rho \bar{\psi}_{2}^{2}-%
\frac{\bar{T}}{u^{2}}\bar{\psi}_{1}^{2})+\frac{4i\sigma _{\bar{z}}}{u}\bar{%
\psi}_{1}\bar{\psi}_{2},\eqno(43a)
$$
$$
G_{2}=4\bar{\rho}(\psi _{1}^{2}-\frac{4\bar{T}}{u^{2}}\psi _{2}^{2}-\frac{4i%
\bar{\sigma}_{\bar{z}}}{u}\psi _{1}\psi _{2}-4\rho \bar{\psi}_{2}^{2}+\frac{%
\bar{T}}{u^{2}}\bar{\psi}_{1}^{2})-\frac{4i\sigma _{\bar{z}}}{u}\bar{\psi}%
_{1}\bar{\psi}_{2},\eqno(43b)
$$
$$
G_{3}=2[4\rho (\psi _{1}\bar{\psi}_{2}+\frac{\bar{T}}{u^{2}}\bar{\psi}%
_{1}\psi _{2})+\frac{2i\sigma _{\bar{z}}}{u}(|\psi _{2}|^{2}-|\psi
_{1}|^{2})],\eqno(43c)
$$
which enter into the following exact forms 
$$
x_{1}(z,\bar{z})=i\int_{\gamma }F_{1}(z^{\prime },\bar{z}^{\prime
})dz^{\prime }+G_{1}(z^{\prime },\bar{z}^{\prime })d\bar{z}^{\prime },%
\eqno(44a)
$$
$$
x_{2}(z,\bar{z})=\int_{\gamma }F_{2}(z^{\prime },\bar{z}^{\prime
})dz^{\prime }+G_{2}(z^{\prime },\bar{z}^{\prime })d\bar{z}^{\prime },%
\eqno(44b)
$$
$$
x_{3}(z,\bar{z})=\int_{\gamma }F_{3}(z^{\prime },\bar{z}^{\prime
})dz^{\prime }+G_{3}(z^{\prime },\bar{z}^{\prime })d\bar{z}^{\prime },%
\eqno(44c)
$$
integrated on an arbitrary curve $\gamma $ in the complex plane, connecting
an initial fixed point to a final one $z$. The conservation laws (41)
guarantee that the quantities $x_{i}$ do not depend on the choice of the
curve $\gamma $ in the complex plane $C$ . Since the functions $x_{i}(z,\bar{%
z})$ are real, they can be considered as components of the radial vector 
$$
{\bf r}(z,\bar{z})=(x_{1},x_{1},x_{1})\eqno(45)
$$
of an orientable, connected surface immersed in $R^{3}$ and locally
parametrized by $z.$ The properties of such surfaces have to be determined,
but they are closely related to the Weingarten surfaces of Class 2. The
reason for  this is based on the observation that the tangent vectors can be
represented by the formulae (32), after the identification $S\rightarrow 
{\bf S}$ , where ${\bf S}$ is a unimodular real vector of ${\bf R}^{3}.$
Than interpreting such a vector as the normal vector to the surface,
expressions (32) are similar to the equations (4.8) in [11], which are the
so - called Lelieuvre formulas for the Weingarten surfaces of Class 2. The
only difference resides in the last terms of both expressions, which in our
case are proportional to the normal vector and not to its derivatives.
However the coefficients also change, while in the Schief case one has the
same coefficient .

\section{The Weierstrass representation for surfaces immersed into $R^{8}$%
}

Quite recently in [9], a generalization of the WR of equations corresponding
to $CP^{2}$ harmonic maps was presented. This generalization allows to study
2-dimensional surfaces immersed in a flat space $R^{8}$ with Euclidean
metric. Analogously, the generalized spin model (31) also gives a new class
of surfaces immersed into $R^{8}.$ To do so, first we rewrite the equation
(31) in terms of the projector as 
$$
(\rho PP_{y})_{y}+(\rho PP_{x})_{x}+\sigma _{x}P_{x}+\sigma _{y}P_{y}=0%
\eqno(46a) 
$$
$$
\rho _{xx}+\rho _{yy}=0,\qquad \sigma _{xx}+\sigma _{yy}=0\eqno(46b) 
$$
where 
$$
P=\frac{1}{1+|\omega |^{2}}\left( 
\begin{array}{ll}
1 & \bar{\omega} \\ 
\omega & |\omega |^{2}
\end{array}
\right) .\eqno(47) 
$$

Now let us consider the equation (46a) when the projector $P$ is a $3\times
3 $ matrix with the structure 
$$
P=\frac{1}{1+|\omega _{1}|^{2}+|\omega _{2}|^{2}}\left( 
\begin{array}{lll}
1 & \bar{\omega}_{1} & \bar{\omega}_{2} \\ 
\omega _{1}{}^{{}} & |\omega _{1}|^{2} & \omega _{1}\bar{\omega}_{2} \\ 
\omega _{2} & \bar{\omega}_{1}\omega _{2} & |\omega _{2}|^{2}
\end{array}
\right) .\eqno(48) 
$$

As in the previously, in this case we have the conservation laws 
$$
K_{z}+M_{\bar{z}}=0,\eqno(49)
$$
with the matrices 
$$
K=[P_{\bar{z}},P]+\sigma _{\bar{z}}P,\qquad M=[P_{z},P]+\sigma _{z}P.%
\eqno(50)
$$
Then, as in section 3, we can introduce two pairs of complex functions $\psi
_{i},\phi _{i}$ as 
$$
\omega _{i}=\frac{\psi _{i}}{\phi _{i}}\eqno(51)
$$
and express $K,M$ in terms of $\psi _{i},\phi _{i}$. In this way, the
equation (46) allows us to define eight new conservation laws. Since the
results of the calculations are the same as in Ref.[9], we will not present
them here. Only we remark that we can contruct eight real quantities $%
X\_\{i\}(z,$ $z)$ which can be used to study various geometrical aspects of
the surface.

\section{Summary and concluding remarks}

The main aim of this paper is to propose a generalization of the
Konopelchenko system (3), involved in the study of constant mean curvature
surfaces. Supplying the general structure (7) with constraints on the
functions $h_{1,}\quad h_{2},\quad \rho $ and $g$ we can recover known
examples and, in particular we obtained the Weierstrass - like
representation (41) - (42) - (43) for surfaces of Weingarten type. However,
the integrability properties of the system (27), its explicit solutions and
the concrete use of formulae (44) is far to be achieved. Finally, we have
shown that this approach can be pursued also in the case of certain surfaces
immersed in $R^{8}.$

\section{Acknowledgments}

The authors are grateful to A.M.Grundland, B.G.Konopelchenko and G.Landolfi
for very helpful discussions. This work was supported in part by MURST of
Italy, INFN-Sezione di Lecce, and INTAS, grant 99-1782. One of the authors
(R.M.) thanks the Department of Physics of the Lecce University for its kind
hospitality.

\end{document}